%
%
%
\input amstex
\documentstyle{amsppt}  
\pageheight{9truein}\pagewidth{6.5truein}
\NoBlackBoxes

%
%
\PSAMSFonts

%
%
\refstyle{A}

%
%
%
\newif\ifxxx
%
\xxxtrue    
%
%
%
%
%
%

\newif\ifian 

\newif\iftextures
\ifian%
\texturestrue%
\else
\fi

\def\boxit#1#2#3{\hbox{\vrule%
\vtop{%
\vbox{\hrule\kern#1%
\hbox{\kern#1$\scriptstyle #2,#3$\kern#1}}%
\kern#1\hrule}%
\vrule}}
\def\ssboxit#1#2#3{\hbox{\vrule%
\vtop{%
\vbox{\hrule\kern#1%
\hbox{\kern#1$\scriptscriptstyle #2,#3$\kern#1}}%
\kern#1\hrule}%
\vrule}}

\def\psm{{\text{\raise 1.5pt\hbox{\kern 0.5pt \vrule height .33pt depth 0pt width 4pt\kern 0.5pt}}}}
\newread\epsffilein    
\newif\ifepsffileok    
\newif\ifepsfbbfound   
\newif\ifepsfverbose   
\newdimen\epsfxsize    
\newdimen\epsfysize    
\newdimen\epsftsize    
\newdimen\epsfrsize    
\newdimen\epsftmp      
\newdimen\pspoints     
\pspoints=1bp          
\epsfxsize=0pt         
\epsfysize=0pt         
\def\epsfbox#1{\global\def\epsfllx{72}\global\def\epsflly{72}%
   \global\def\epsfurx{540}\global\def\epsfury{720}%
   \def\lbracket{[}\def\testit{#1}\ifx\testit\lbracket
   \let\next=\epsfgetlitbb\else\let\next=\epsfnormal\fi\next{#1}}%
\def\epsfgetlitbb#1#2 #3 #4 #5]#6{\epsfgrab #2 #3 #4 #5 .\\%
   \epsfsetgraph{#6}}%
\def\epsfnormal#1{\epsfgetbb{#1}\epsfsetgraph{#1}}%
\def\epsfgetbb#1{%
%
%
\openin\epsffilein=#1
\ifeof\epsffilein\errmessage{I couldn't open #1, will ignore it}\else
%
%
   {\epsffileoktrue \chardef\other=12
    \def\do##1{\catcode`##1=\other}\dospecials \catcode`\ =10
    \loop
       \read\epsffilein to \epsffileline
       \ifeof\epsffilein\epsffileokfalse\else
%
%
          \expandafter\epsfaux\epsffileline:. \\%
       \fi
   \ifepsffileok\repeat
   \ifepsfbbfound\else
    \ifepsfverbose\message{No bounding box comment in #1; using defaults}\fi\fi
   }\closein\epsffilein\fi}%
%
%
\def\epsfsetgraph#1{%
   \epsfrsize=\epsfury\pspoints
   \advance\epsfrsize by-\epsflly\pspoints
   \epsftsize=\epsfurx\pspoints
   \advance\epsftsize by-\epsfllx\pspoints
%
%
   \epsfsize\epsftsize\epsfrsize
   \ifnum\epsfxsize=0 \ifnum\epsfysize=0
      \epsfxsize=\epsftsize \epsfysize=\epsfrsize
%
%
     \else\epsftmp=\epsftsize \divide\epsftmp\epsfrsize
       \epsfxsize=\epsfysize \multiply\epsfxsize\epsftmp
       \multiply\epsftmp\epsfrsize \advance\epsftsize-\epsftmp
       \epsftmp=\epsfysize
       \loop \advance\epsftsize\epsftsize \divide\epsftmp 2
       \ifnum\epsftmp>0
          \ifnum\epsftsize<\epsfrsize\else
             \advance\epsftsize-\epsfrsize \advance\epsfxsize\epsftmp \fi
       \repeat
     \fi
   \else\epsftmp=\epsfrsize \divide\epsftmp\epsftsize
     \epsfysize=\epsfxsize \multiply\epsfysize\epsftmp   
     \multiply\epsftmp\epsftsize \advance\epsfrsize-\epsftmp
     \epsftmp=\epsfxsize
     \loop \advance\epsfrsize\epsfrsize \divide\epsftmp 2
     \ifnum\epsftmp>0
        \ifnum\epsfrsize<\epsftsize\else
           \advance\epsfrsize-\epsftsize \advance\epsfysize\epsftmp \fi
     \repeat     
   \fi
%
%

\iftextures
   \ifepsfverbose\message{#1: width=\the\epsfxsize, height=\the\epsfysize}\fi
   \epsftmp=10\epsfxsize \divide\epsftmp\pspoints
   \vbox to\epsfysize{\vfil\hbox to\epsfxsize{%
      \special{illustration #1 scaled \number\epsfscale}
      \hfil}}%
\else
   \ifepsfverbose\message{#1: width=\the\epsfxsize, height=\the\epsfysize}\fi
   \epsftmp=10\epsfxsize \divide\epsftmp\pspoints
   \vbox to\epsfysize{\vfil\hbox to\epsfxsize{%
      \includegraphics{#1}%
     \hfil}}%
\fi
\epsfxsize=0pt\epsfysize=0pt\epsfscale=1000 }%

%
%
{\catcode`\%=12 \global\let\epsfpercent=
%
%
\long\def\epsfaux#1#2:#3\\{\ifx#1\epsfpercent
   \def\testit{#2}\ifx\testit\epsfbblit
      \epsfgrab #3 . . . \\%
      \epsffileokfalse
      \global\epsfbbfoundtrue
   \fi\else\ifx#1\par\else\epsffileokfalse\fi\fi}%
%
%
\def\epsfgrab #1 #2 #3 #4 #5\\{%
   \global\def\epsfllx{#1}\ifx\epsfllx\empty
      \epsfgrab #2 #3 #4 #5 .\\\else
   \global\def\epsflly{#2}%
   \global\def\epsfurx{#3}\global\def\epsfury{#4}\fi}%
%
%
%
%

\newcount\epsfscale    
\newdimen\epsftmpp     
\newdimen\epsftmppp    
\newdimen\epsfM        
\newdimen\sppoints     
\epsfscale=1000        
\sppoints=1000sp       
\epsfM=1000\sppoints
%
\def\computescale#1#2{%
  \epsftmpp=#1 \epsftmppp=#2
  \epsftmp=\epsftmpp \divide\epsftmp\epsftmppp  
  \epsfscale=\epsfM \multiply\epsfscale\epsftmp 
  \multiply\epsftmp\epsftmppp                   
  \advance\epsftmpp-\epsftmp                    
  \epsftmp=\epsfM                               
  \loop \advance\epsftmpp\epsftmpp              
    \divide\epsftmp 2                           
    \ifnum\epsftmp>0
      \ifnum\epsftmpp<\epsftmppp\else           
        \advance\epsftmpp-\epsftmppp            
        \advance\epsfscale\epsftmp \fi          
  \repeat
  \divide\epsfscale\sppoints}
\def\epsfsize#1#2{%
  \ifnum\epsfscale=1000
    \ifnum\epsfxsize=0
      \ifnum\epsfysize=0
      \else \computescale{\epsfysize}{#2}
      \fi
    \else \computescale{\epsfxsize}{#1}
    \fi
  \else
    \epsfxsize=#1
    \divide\epsfxsize by 1000 \multiply\epsfxsize by \epsfscale
  \fi}




\def\drawsmash#1{\relax
\ifmmode%
\typeout{Warning:Macro drawsmash used in math mode}%
\fi%
\setbox0=\hbox{#1}\ht0=0pt\dp0=0pt\box0}

\newcount\figcount
\figcount=0
\newbox\drawing
\newcount\drawbp
\newdimen\drawx
\newdimen\drawy
\newdimen\ngap
\newdimen\sgap
\newdimen\wgap
\newdimen\egap

\def\drawbox#1#2#3{\vbox{
  \setbox\drawing=\vbox{\offinterlineskip\epsfbox{#2.eps}\kern 0pt}
  \drawbp=\epsfurx
  \advance\drawbp by-\epsfllx\relax
  \multiply\drawbp by #1
  \divide\drawbp by 100
  \drawx=\drawbp truebp
  \ifdim\drawx>\hsize\drawx=\hsize\fi
  \epsfxsize=\drawx
  \setbox\drawing=\vbox{\offinterlineskip\epsfbox{#2.eps}\kern 0pt}
  \drawx=\wd\drawing
  \drawy=\ht\drawing
  \ngap=0pt \sgap=0pt \wgap=0pt \egap=0pt
  \setbox0=\vbox{\offinterlineskip
    \box\drawing \ifgridlines\drawgrid\drawx\drawy\fi #3}
  \kern\ngap\hbox{\kern\wgap\box0\kern\egap}\kern\sgap}}

\def\draw#1#2#3{
\setbox\drawing=\drawbox{#1}{#2}{#3}
\vskip 6pt \centerline{\ifgridlines\boxgrid\drawing\fi\box\drawing}\vskip 6pt%
}

\newif\ifgridlines
\newbox\figtbox
\newbox\figgbox
\newdimen\figtx
\newdimen\figty

\newdimen\bwd
\bwd=2sp 

\def\hhline#1{\vbox{\drawsmash{\hbox to #1{\leaders\hrule height \bwd\hfil}}}}
\def\vline#1{\hbox to 0pt{%
  \hss\vbox to #1{\leaders\vrule width \bwd\vfil}\hss}}

\def\clap#1{\hbox to 0pt{\hss#1\hss}}
\def\vclap#1{\vbox to 0pt{\offinterlineskip\vss#1\vss}}

\def\hstutter#1#2{\hbox{%
  \setbox0=\hbox{#1}%
  \hbox to #2\wd0{\leaders\box0\hfil}}}

\def\vstutter#1#2{\vbox{
  \setbox0=\vbox{\offinterlineskip #1}
  \dp0=0pt
  \vbox to #2\ht0{\leaders\box0\vfil}}}

\def\crosshairs#1#2{
  \dimen1=.002\drawx
  \dimen2=.002\drawy
  \ifdim\dimen1<\dimen2\dimen3\dimen1\else\dimen3\dimen2\fi
  \setbox1=\vclap{\vline{2\dimen3}}
  \setbox2=\clap{\hhline{2\dimen3}}
  \setbox3=\hstutter{\kern\dimen1\box1}{4}
  \setbox4=\vstutter{\kern\dimen2\box2}{4}
  \setbox1=\vclap{\vline{4\dimen3}}
  \setbox2=\clap{\hhline{4\dimen3}}
  \setbox5=\clap{\copy1\hstutter{\box3\kern\dimen1\box1}{6}}
  \setbox6=\vclap{\copy2\vstutter{\box4\kern\dimen2\box2}{6}}
  \setbox1=\vbox{\offinterlineskip\box5\box6}
  \drawsmash{\vbox to #2{\hbox to #1{\hss\box1}\vss}}}

\def\boxgrid#1{\rlap{\vbox{\offinterlineskip
  \setbox0=\hhline{\wd#1}
  \setbox1=\vline{\ht#1}
  \drawsmash{\vbox to \ht#1{\offinterlineskip\copy0\vfil\box0}}
  \drawsmash{\vbox{\hbox to \wd#1{\copy1\hfil\box1}}}}}}

\def\drawgrid#1#2{\vbox{\offinterlineskip
  \dimen0=\drawx
  \dimen1=\drawy
  \divide\dimen0 by 10
  \divide\dimen1 by 10
  \setbox0=\hhline\drawx
  \setbox1=\vline\drawy
  \drawsmash{\vbox{\offinterlineskip
    \copy0\vstutter{\kern\dimen1\box0}{10}}}
  \drawsmash{\hbox{\copy1\hstutter{\kern\dimen0\box1}{10}}}}}

\def\figtext#1#2#3#4#5{
  \setbox\figtbox=\hbox{#5}
  \dp\figtbox=0pt
  \figtx=-#3\wd\figtbox \figty=-#4\ht\figtbox
  \advance\figtx by #1\drawx \advance\figty by #2\drawy
  \dimen0=\figtx \advance\dimen0 by\wd\figtbox \advance\dimen0 by-\drawx
  \ifdim\dimen0>\egap\global\egap=\dimen0\fi
  \dimen0=\figty \advance\dimen0 by\ht\figtbox \advance\dimen0 by-\drawy
  \ifdim\dimen0>\ngap\global\ngap=\dimen0\fi
  \dimen0=-\figtx
  \ifdim\dimen0>\wgap\global\wgap=\dimen0\fi
  \dimen0=-\figty
  \ifdim\dimen0>\sgap\global\sgap=\dimen0\fi
  \drawsmash{\rlap{\vbox{\offinterlineskip
    \hbox{\hbox to \figtx{}\ifgridlines\boxgrid\figtbox\fi\box\figtbox}
    \vbox to \figty{}
    \ifgridlines\crosshairs{#1\drawx}{#2\drawy}\fi
    \kern 0pt}}}}


\def\hpad#1#2#3{\hbox{\kern #1\hbox{#3}\kern #2}}
\def\vpad#1#2#3{\setbox0=\hbox{#3}\dp0=0pt\vbox{\kern #1\box0\kern #2}}



\def\stack#1#2#3{\vbox{\offinterlineskip
  \setbox2=\hbox{#2}
  \setbox3=\hbox{#3}
  \dimen0=\ifdim\wd2>\wd3\wd2\else\wd3\fi
  \hbox to \dimen0{\hss\box2\hss}
  \kern #1
  \hbox to \dimen0{\hss\box3\hss}}}


\def\hexp#1{%
  \setbox0=\hbox{${}^{#1}$}%
  \hbox to .5\wd0{\box0\hss}}

\define\sltwo{\operatorname {SL}_2}
\define\af #1.{\Bbb A^{#1}}
\define\au#1.{\operatorname {Aut}\,(#1)}
\define\pr #1.{\Bbb P^{#1}}
\define\pic#1.{\operatorname {Pic}\,(#1)}
\define\bbbq{{\Bbb Q}}
\define\bbbv{{\Bbb V}}
\define\bbbr{{\Bbb R}}

\define\qpic#1.{{\pic #1.}_{\bbbq}}
\define\ol#1{\overline{#1}}
\define\ses#1.#2.#3.{0\longrightarrow #1 \longrightarrow #2 \longrightarrow #3
\longrightarrow 0}

\define\mg{\overline { M}_g}
\define\otn{\{1,2,\dots,n\}}
\define\vmgn#1.#2.{\overline{M}_{#1,#2}}
\define\vmg #1.{\overline{M}_{#1}}
\define\mgn{\overline{M}_{g,n}}
\define\mon{\vmgn 0. n.}
\define\vln#1.{\overline{L}_{#1}}
\define\rgn{\overline{R}_{g,n}}
\define\vrgn#1.#2.{\overline{R}_{#1,#2}}
\define\fgn{\overline{F}_{g,n}}
\define\bmgn#1.{\partial{\mgn #1.}}

\define\ex#1.{{\Bbb E}(#1)}
\define\mc#1.{\overline{NE}_1(#1)}
\define\md#1.{\overline{NE}^1(#1)}
\define\tmc#1.{\overline{NE}_1(#1)}
\define\tpic#1.{\operatorname{Pic}(#1)}
\define\tqpic#1.{{\tpic #1.}_{\bbbq}}
\define\omc#1.{{NE}_1(#1)}
\define\tomc#1.{{NE}_1(#1)}
\define\rmc#1.#2.{\overline{NE}_1(#1/#2)}

\define\sym#1.#2.{\operatorname {Sym}^{#1}(#2)}
\define\rpic#1.#2.#3.{\operatorname {Pic}_{#1/#2}^{#3}}


\define\irr{\operatorname{irr}}
\define\birr{b_{\irr}}
\define\dirr{\delta_{\irr}}
\define\ct{{\Cal T}}

\define\cz{{\Cal Z}}
\define\te{\tilde{E}}

%
%
%
%

%
%
%
%

\topmatter
\title Towards the ample cone of $\overline{M}_{g,n}$ \endtitle
\author Angela Gibney, Sean Keel and Ian Morrison \endauthor

\leftheadtext{Angela Gibney, Sean Keel and Ian Morrison}%
\rightheadtext{Towards the ample cone of $\overline{M}_{g,n}$}%

\address Department of Mathematics, University of Michigan, 
Ann Arbor, MI 48109 \endaddress

\email agibney\@math.lsa.umich.edu \endemail

\address Department of Mathematics, University of Texas at Austin, 
Austin, TX 78712 \endaddress

\email keel\@fireant.ma.utexas.edu \endemail

\address Department of Mathematics, Fordham University,
Bronx, NY 10458 \endaddress

\email morrison\@fordham.edu \endemail


\abstract
  In this paper we study the ample cone of the moduli space $\mgn$
of stable $n$-pointed curves of genus $g$. 
Our motivating
conjecture is that a divisor on $\mgn$ is ample iff it has positive
intersection with all $1$-dimensional strata 
(the components of the locus of
curves with at least $3g+n-2$ nodes). 
This translates into
a simple conjectural
description of the cone by linear inequalities, and, 
as all the $1$-strata are rational, includes
the conjecture that the Mori cone is polyhedral and 
generated by rational curves.
Our main result is that the conjecture holds iff it
holds for $g=0$. More precisely, there is a natural finite
map $r: \vmgn 0. 2g+n. \rightarrow \mgn$ whose image
is the locus $\rgn$ of curves with all components rational.
Any $1$-strata either lies in $\rgn$ or is numerically 
equivalent to a family $E$ of elliptic tails and we show that a 
divisor $D$ is nef iff $D \cdot E \geq  0$ and $r^*(D)$ is 
nef. We also give results on
contractions (i.e. morphisms with connected fibers to
projective varieties) of $\mgn$ for $g \geq 1$ showing
that any fibration factors through a tautological
one (given by forgetting points) and that the
exceptional locus of any birational contraction
is contained in the boundary.

\endabstract

\thanks During this research, the first two authors received partial
support from a Big XII faculty research grant, and a 
grant from the Texas Higher Education Coordinating Board.\endgraf
The first author received further
partial support from the Clay Math Institute, and
the second from the NSF.\endgraf
The third author's research was
partially supported by a Fordham University Faculty Fellowship and
by grants from the Centre de
Recerca Matem\'atica for a stay in
Barcelona and from the Consiglio Nazionale di Ricerche for stays in Pisa
and Genova. 
\endthanks
 

\dedicatory To Bill Fulton on his sixtieth birthday \enddedicatory


\endtopmatter

%
%
%
%

%
%
%
%

\document

\subhead \S 0 Introduction and Statement of Results \endsubhead

The moduli space of stable curves is among the most studied
objects in algebraic geometry. Nonetheless, its birational
geometry remains largely a mystery and most Mori theoretic
problems in the area are entirely open. Here we consider
one of the most basic:

\proclaim{0.1 Question (Mumford)} What are the ample divisors
on $\mgn$?
\endproclaim

There is a stratification of $\mgn$ by
topological type where the codimension $k$ strata are the
irreducible components of the locus parameterizing
pointed curves with at least $k$ singular points. An ample
divisor must intersect any one dimensional stratum positively.
It is thus natural to consider the following conjecture.

\proclaim{0.2 Conjecture: $F_1(\mgn)$} A divisor
on $\mgn$ is ample iff it has positive intersection
with all $1$-dimensional strata. In other words, any
effective curve in $\mgn$ is numerically equivalent to
an effective combination of $1$-strata. Said differently still,
 every extremal ray of the Mori
cone of effective curves $\mc {\mgn}.$ is generated
by a one dimensional stratum.
\endproclaim

$F_1(\vmgn 0.n.)$ was previously conjectured by
Fulton, whence our choice of notation. Special
cases of $F_1(\mon)$ are known (see \cite{KeelMcKernan96},
\S 7) but
the general case remains open. The central observation of this work is that the case of higher
genus is no harder. Indeed, $\mc {\mgn}.$ is naturally
a quotient of $\mc {\vmgn 0. 2g+n.}. \times {\Bbb R}^{\geq 0}$
and Conjecture (0.2) is true for all 
$g$ iff it is true for $g=0$: see (0.3 -0.5).
Further, the second formulation of conjecture (0.2) holds except
possibly for very 
degenerate families in $\mgn$: see (0.6). 
Moreover, we are able to give strong results on the 
contractions of $\mgn$
(i.e. morphisms with geometrically connected fibers from it to 
other projective varieties).  For example we show that
for $g \geq 2$ the only fibrations of $\mgn$ are compositions
of a tautological fibration, given by dropping
points, with a birational morphism, see (0.9-0.11).
Here are precise statements of our results:

We work over an algebraically closed field
of any characteristic other than $2$ (this last assumption is used
only at the end of the proof of (3.5)).
 
By the locus of {\it flag
curves}, we shall mean the image $\fgn$ of the 
morphism $f: \vmgn 0. g + n./S_g \rightarrow \mgn$
induced by gluing on $g$ copies of the pointed
rational elliptic curve at $g$ points (which the symmetric
group $S_g \subset S_{g+n}$ permutes). The map $f$ is the normalization
of $\fgn$.

\proclaim{0.3 Theorem} $g \geq 1$. 
A divisor $D \in \tpic \mgn.$ is
nef iff $D$ has non-negative intersection
with all the one dimensional strata and $D|_{\fgn}$ is nef.
In particular, $F_1(\vmgn 0. g+n./S_g)$ implies $F_1(\mgn)$.
\endproclaim

Using the results of  \cite{KeelMcKernan96}, Theorem (0.3) has the following consequence.

\proclaim{0.4 Corollary}
$F_1(\mgn)$ holds in characteristic zero
as long as $g+n \leq 7$ and when $n=0$ for
$g \leq 11$.
\endproclaim

For  $g \leq 4$, $F_1(\mg)$ was obtained previously by Faber.

We call $E \subset \mgn$ ($g \geq 1$) a 
{\it family of elliptic tails} if it
is the image of the map
$\vmgn 1.1. \rightarrow \mgn$ obtained by attaching
a fixed $n+1$-pointed curve of genus $g-1$ to the moving 
pointed elliptic curve.
Any two families of elliptic tails are numerically equivalent.
In (2.2), we will see that all the $1$-strata except for $E$ are
numerically equivalent to families of rational curves. Note that 
 we abuse language and refer to a curve
as rational if all of its irreducible components are rational. 
Thus it
is immediate from (0.3) that the Mori cone is generated by
$E$ together with curves in $\rgn \subset \mgn$, the locus of rational
curves. The locus $\rgn$ is itself (up to normalization) a quotient of
$\vmgn 0. 2g + n.$: By deformation theory 
$\rgn$ is the closure of the locus of
irreducible $g$-nodal (necessarily rational) curves. The
normalization of such a curve is a
smooth rational curve with $2g+n$ marked points.
Thus, if
$G \subset S_{2g}$ is the subgroup of permutations commuting with the
product of $g$ transpositions
$(12)(34)\dots(2g-1\,2g)$, the normalization of $\rgn$ is naturally
identified with
$\vmgn 0. 2g+n./G$. Thus (0.3) implies the following result:

\proclaim{0.5 Corollary} $D \in \tpic \mgn.$ is nef 
iff $D|_{\rgn \cup E}$ is nef. Equivalently,
the natural map
$$
\mc { {\vmgn 0. 2g +n./G} \cup E}. = {\mc {\vmgn 0. 2g+n./G}.} \times
{\bbbr}^{\geq 0} \rightarrow \tmc {\mgn}.
$$
is surjective.
\endproclaim

In fact we prove the following strengthening of (0.3):

\proclaim{0.6 Theorem} Let $g \geq 1$ and let $N \subset \tmc \mgn.$ be
the subcone generated by the strata (1-5) of (2.2) for $g \ge 3$, by  (1)
and (3-5) for $g=2$, and by (1) and (5) for $g=1$. Then
$N$ is the 
subcone generated by curves $C \subset \mgn$ whose
associated family of curves has no moving smooth rational
components. Hence,
$$
\tomc {\mgn}. = N + \mc {\vmgn 0. g+n./S_g}.
$$
\endproclaim

(See (1.2) for the meaning of {\it moving component}.) An analogous result
for families whose general member has at most one 
singularity is given in  \cite{Moriwaki00,A}. 

Of course Theorem (0.6) implies (the second formulation) of
Conjecture (0.2) for any family $C \subset \mgn$ with
no moving smooth rational component. In fact, 
 $F_1(\vmgn 0. 7.)$ holds, and so Conjecture (0.2) holds as long as 
there is no moving smooth rational component 
containing at least $8$ distinguished (i.e. 
marked or singular) points. 

Theorem (0.3) also has a converse which we state as follows.

\proclaim{0.7 Theorem} Every nef line bundle on
$\vmgn 0. g+n./S_g$ is the pullback of a nef line bundle
on $\mgn$ and  $F_1(\vmgn 0. g+n./S_g)$ is equivalent
to $F_1(\mgn)$. In particular, $F_1(\mg)$ is
equivalent to $F_1(\vmgn 0. g. /S_g)$ and $F_1(\vmgn 1. n.)$
is equivalent to $F_1(\vmgn 0. n+1.)$. 
\endproclaim

\proclaim{0.8 Theorem} The Mori cone of $\fgn$ is
a face of the Mori cone of $\mgn$:
there is a nef divisor $D$ such that
$$
D^{\perp} \cap \mc{\mgn}. = \mc {\vmgn 0. g+n.}/S_g. .
$$
\endproclaim

\remark{Remark} An interesting problem is to determine
whether or not 
$F_1(\mg)$ for $g \geq 1$
implies $F_1(\mon)$ for all $n \geq 3$. For example, if one maps
$\mon$ into $\mg$ by gluing one pointed curves of very
different genera, then the pullback on the Picard groups will
be surjective. However, it is not clear that this will be true for the
nef cones.
\endremark

We note that in a number of respects using $\fgn$ to reduce the general conjecture
is more desirable than using $\rgn$. Analogues of Theorems 
(0.7) and (0.8) fail for $\rgn$, indicating
that $\fgn$ is the more natural locus from
the Mori theoretic point of view (e.g. one is led to
consider $\fgn$ by purely combinatorial consideration 
of the Mori cone). Further 
the birational geometry of $\fgn$ is much simpler than
that of $\rgn$. For example, when $n=0$,
the Picard group of $\vmgn 0. g./S_g$ has rank roughly
$g/2$ and is freely generated
by the boundary divisors, while
$\pic {\vmgn 0. 2g./G}.$ has rank roughly $g^2/2$
and there is a relation among the 
boundary classes. Of greater geometric significance is
the contrast between the
cones of effective divisors:
$\md {\vmgn 0. g./S_g}.$ is simplicial and is generated
by the boundary divisors (see \cite{KeelMcKernan96}),
while the structure of $\md {\vmgn 0. 2g.}/G .$ is
unclear and it is definitely not generated by boundary divisors
(a counter-example is noted below).

It is straightforward to enumerate the possibilities
for the one dimensional strata, and then express (0.2)
as a conjectural description of the ample
cone as an intersection of explicit half spaces. 
This description  is given for $n=0$ 
in \cite{Faber96} and for $g=0$
in \cite{KeelMcKernan96}. We treat the general case in \S 2.

A fundamental problem in birational geometry is to study
morphisms of a given variety to other varieties. In
the projective category any such morphism is given
by a semi-ample divisor -- i.e. a divisor such that the linear
system of some positive multiple is basepoint free. A
semi-ample divisor is necessarily nef. This implication
is one of the main reasons for considering nef
divisors. Correspondingly, one of the main reasons that (0.1) is 
interesting is its connection to:

\proclaim{Question} What are all the contractions (i.e. morphisms
with geometrically connected fibers to projective varieties) 
of $\mgn$?
\endproclaim

We have a number of results in this direction.
Recall a definition from \cite{Keel99}: 
For a nef divisor $D$ on a variety
$X$ a subvariety $Z \subset X$ is called $D$-exceptional
if the $D$-degree of $Z$ is zero, or equivalently, $D|_Z$
is not big. The exceptional locus of $D$, denoted
$\ex D.$, is the union of the exceptional subvarieties.
If $D$ is semi-ample
and big, this is the exceptional 
locus for the associated birational map.
  
\proclaim{0.9 Theorem}
Let $D \in \tpic{\mgn}.$ be a nef divisor,
$g \geq 2$ (resp. $g=1$). Either
$D$ is the pull back of a nef divisor 
on $\vmgn g. n-1.$ via one of the tautological
projections (resp. is the tensor product 
of pullbacks of nef divisors on $\vmgn 1.S.$
and $\vmgn 1.S^c.$ via the tautological projection
for some subset $S \subset N$)
or $D$ is big and 
$\ex D. \subset \partial \mgn$.
\endproclaim

\proclaim{0.10 Corollary} For $g \geq 2$ any
fibration of $\mgn$ (to a projective variety)
factors through a projection to some
$\vmgn g. i.$ ($i < n$), while  $\mg$ has no non-trivial fibrations.
\endproclaim

The above corollary and its proof are part of the first author's
Ph.D. thesis, \cite{Gibney00}, written under the direction of
the second author.

\proclaim{0.11  Corollary} $g \geq 1$.
Let $f: \mgn \rightarrow X$ be a birational
morphism to a projective variety.
Then the exceptional
locus of $f$ is contained in $\partial \mgn$.
In particular $X$ is again a compactification of $M_{g,n}$.
\endproclaim

Note that our results give some support to 
the conjecture that the nef cone of $\mgn$ behaves
as if it were described by Mori's cone theorem (i.e. like the
cone of a log Fano variety) with dual cone generated by finitely
many rational curves and having every face contractible. This is 
surprising, since for example $\mgn$ is usually 
(e.g. for $g \ge 24$) of general
type, and conjecturally (by Lang, see \cite{CHM97}) 
{\it no} rational curve meets the interior of
$\mgn$ for $g \geq 2$ and $n$ sufficiently big.
Of course, our main results reduce Conjecture (0.2) to genus
$0$, where from initial Mori theoretic considerations it is
much more plausible. 

We note that prior to this work very few nef (but not ample)
line bundles on $\mgn$ were known.
Myriad examples can be constructed using (0.3) as
we indicate in \S 6.

The above results indicate that (regular) contractions give only
a narrow view of the geometry of $\mgn$.  For example, (0.10)
gives that the only fibrations are 
(induced from) the tautological ones, and by (0.11), birational
morphisms only affect the boundary. On the other hand the
birational geometry should be very rich. For example, as is seen in (6.4), all but one
of the elementary
(i.e. relative Picard number one) extremal contractions
of $\mg$ 
are small (i.e. isomorphisms in codimension one), 
so there should be a wealth of flips. Also one expects
in many cases to find interesting rational fibrations.   For example,
when $g+1$ has more than one factorization, the Brill-Noether
divisor should give a positive dimensional linear system on
$\mg$ which is conjecturally (see \cite{HarrisMorrison98,6.63})
not big. 

One consequence of (0.9) and the fact
that $\ex {\lambda}. = \partial \mg$ is that $\partial \mg$ is intrinsic.

\proclaim{0.12 Corollary}  Any automorphism of $\mg$ must preserve the boundary.
 \endproclaim

In view of our results, $F_1(\mon)$ is
obviously of compelling interest. It is natural to wonder:

\proclaim{0.13 Question} If a divisor on $\mon$ has non-negative
intersection with all one dimensional strata, does it follow
that the divisor is linearly equivalent to an effective
combination of boundary divisors?
\endproclaim

A positive answer to (0.13) would imply (0.2). The analogous
property does hold on $\vmgn 1.n.$, (see (3.3)) but the 
question in that case is vastly simpler 
as the boundary divisors are linearly independent.

It is not true that every effective divisor on $\mon$ is
an effective sum of boundary divisors.  In other words,  $F^1(\mon)$
is false.
In particular, the second author has shown
(in support of the slope conjecture, \cite{HarrisMorrison98,6.63})
that the pullback of the Brill-Noether divisor from
$\mg$ to $\vmgn 0. 2g.$ is a counter-example to
$F^1(\vmgn 0. 2g+n./G)$ for any $g \geq 3$
with $g+1$ composite.  In fact, he has shown that for $g=3$ the unique
non-boundary component of this pullback gives a {\it new}
vertex of $\md {\vmgn 0. 6./G}.$. 

Question (0.13) can be
formulated as an elementary combinatorial question of whether
one explicit polyhedral cone 
is contained in another. 
As this restatement
might be of interest to
someone without knowledge of $\mg$, we give this formulation:

\proclaim{0.14 Combinatorial Formulation of (0.13)} Let
$\bbbv$ be the $\bbbq$-vector space that is spanned
by symbols $\delta_T$ for
each subset $T \subset \otn$ subject to the relations
\roster
\item $\delta_T = \delta_{T^c}$ for all $T$;
\item $\delta_T =0$ for $|T| \leq 1$; and
\item for each $4$ element subset $\{i,j,k,l\} \subset \otn$
$$\sum_{i,j \in T,k,l \in T^c} \delta_T = \sum_{i,k \in T,j,l \in T^c}
\delta_{T}\,.$$
\endroster
Let $N \subset \bbbv$ be the set of elements $\sum b_T \delta_T$
satisfying
$$
b_{I\cup J} + b_{I \cup K} + b_{I \cup L} \geq b_I + b_J + b_K +
b_L,
$$
for each partition of $\otn$ into $4$ disjoint subsets
$I,J,K,L$. Let $E \subset \bbbv$ be the affine hull of the
$\delta_T$.

Question: Is $N \subset E$?
\endproclaim

In theory (0.14) can be checked, for a given $n$, by computer.
For $n \leq 6$ we have done so, with considerable
help from Maroung Zou,
using the program Porta. Unfortunately the computational
complexity is enormous, and beyond our machine's capabilities
already for $n=7$. These cases ($n \leq 6$) were proved previously
(by hand) by Faber.  \medskip

The remainder of the paper is organized as follows. In \S 1, we fix 
notation for various 
boundary divisors and gluing maps. In \S 2, we give a 
small set of generators 
for the cone spanned by one-dimensional strata (2.2) and  
inequalities which cut out the dual cone of divisors (2.1) 
--- we term each the Faber cone. The main result of \S 3 is (3.1)
stating that a class nef on the boundary of $\mgn$ is nef. 
In \S 4 the results from the previous section are
used
to deduce (0.6) and hence (0.3).  We also study the nefness and 
exceptional loci
of certain divisors 
in order to deduce (0.7) and (0.8).  The  proof
of (0.9) is contained in \S 5 following preparatory lemmas 
(5.1-5.3) dealing with divisors
on $C^n/G$ where $C$ is a smooth curve with automorphism group $G$. 
In \S 6, we collect some more ad hoc results. In particular, we answer
a question posed by Faber
(6.2) and recover the classical ampleness result of Cornalba-Harris (6.3).
Finally, \S 7 contains a geometric reformulation of $F_1(\mon)$ and a
review of the evidence (to our minds considerable) for (0.2).

\subhead Thanks \endsubhead

We have had interesting and fruitful discussions with various
people on topics related to this paper, including  X.~Chen, 
F.~Cukierman, J.~Harris,  B.~Hassett, J.~Koll\'ar, A.~Logan, J.~McKernan,  
A.~Moriwaki, R.~Pandharipande, 
S.~Popescu, W.~Rulla, S.~L.~Tan, G.~Ziegler, and M.~Zou.
J.~Koll\'ar gave us a counter-example to a theorem claimed
in an earlier version of the paper. 

Especially inspiring and enlightening were lengthy discussions with Carel
Faber. In particular Faber showed us
his proofs of $F_1(\mg)$, $g \leq 4$, and helped us derive
the results of \S 2. 

\subhead \S 1 Notation \endsubhead

For the most part we use 
standard notation for divisors, line bundles, and
loci on $\mgn$. See e.g. \cite{Faber96} and \cite{Faber97}with
the following possible exception.
By $\omega_i$ on $\mgn$, $g \geq 2$, we
mean the pullback of the relative dualizing sheaf
of the universal curve over $\mg$ 
by the projection given by dropping all but
the $i^{th}$ point. Note this is not (for $n \geq 3$)
the relative
dualizing sheaf for the map to $\vmgn g. {n-1}.$
given by dropping the $i^{th}$ point (the symbol
is used variously in the literature).

We note that the $\bbbq$-Picard group of $\mgn$
is the same in all characteristics, by \cite{Moriwaki01}.

We let $N := \otn$.
The cardinality of a finite set $S$ is denoted $|S|$.

For a divisor $D$ on a variety $X$ and $Y$ a closed subset
of $X$, we
say that $D$ is {\it nef outside} $Y$,  if $D \cdot C \geq 0$ for all irreducible
curves $C \not \subset Y$.

To obtain a symmetric description of boundary divisors  on $\mon$, we use
 partitions $\ct = [T,T^c]$ of
$N$ into disjoint subsets, each with at least two
elements, rather than subsets $T$. We write the corresponding
boundary divisor as $\delta_{0,\ct}$, or $\delta_{\ct}$.
For a partition $Q$ of a
subset
$S \subset N$ we write $\ct > Q$ provided the
equivalence relation induced by $Q$ on $S$ refines that obtained by
restricting
$\ct$ to
$S$. We write $Q = \ct|_S$ provided these two equivalence relations
are the same. 

We write $B_k$ for the sum of boundary divisors on $\mon/S_n$ that 
are the image of any $\delta_{\ct}$ with $|T| =k$ (thus $B_k = B_{n-k}$).
We will abuse notation by using the same expression
for its inverse image (with reduced structure) on $\mon$,
$\sum_{\ct, |T|=k} \delta_{\ct}$.

We make repeated use of the standard {\it product decomposition}
for strata as a finite image of products of various $\mgn$. For
precise details see \cite{Keel99,pg. 274} and \cite{Faber97}. 
To describe this decomposition
for boundary divisors we use the notations
$$
\tilde{\Delta}_{i,S}:=
\vmgn i. {S \cup *}. \times \vmgn g-i. {S^c \cup *}. \twoheadrightarrow
\Delta_{i,S}
$$
and
$$
\tilde{\Delta}_{\irr}:=
\vmgn g-1. {N \cup p \cup q}. \twoheadrightarrow \Delta_{\irr}.
$$
We refer to $\vmgn i. S \cup *.$ (and any analogous term
for a higher codimension stratum) as a {\it factor} of the stratum.

The key fact about these maps used here is the

\proclaim{1.1 Lemma} 
The pullback to $\tilde{\Delta}_{i,S}$ of any line bundle 
is numerically equivalent to a tensor product of
unique line bundles from the two factors. The given line bundle is
nef on $\Delta_{i,S}$ iff each of the line bundles on the factors
is nef. Dually, let $C$ be any curve on the product, and
$C_i,C_{g-i}$ be its images on the two factors (with multiplicity
for the pushforward of cycles) which we also view as 
curves in $\mgn$ by the usual device of 
gluing on a fixed curve. Then, $C$ and
$C_i + C_{g-i}$ are numerically equivalent.
\endproclaim
\demo{Proof} The initial statement implies the other statements,
and follows from the explicit formulae of \cite{Faber97}.\qed \enddemo

Any curve $E$ in $\mgn$ induces a decomposition of the curves it
parameterizes into a subcurve fixed in the family $E$ and a moving subcurve.
Arguing inductively, (1.1) yields,

\proclaim{1.2 Corollary}
Every curve in $\mgn$ is numerically
equivalent to an effective combination of curves whose moving subcurves are
all generically irreducible.
\endproclaim

\subhead \S 2. The Faber cone \endsubhead

In this section we consider the subcone of
$\mc {\mgn}.$ generated by one dimensional strata, and
its dual. We refer to either as the {\bf Faber cone}, appending
{\it of curves}, or {\it of divisors} if confusion is possible. We will
call a divisor {\bf $F$-nef} if it lies in the Faber cone. Of course
(0.2) is then the statement that $F$-nef implies nef, or equivalently,
the Mori and Faber cones of $\mgn$ are the same.

The next result describes the Faber cone of divisors as an
intersection of half spaces.
In order to give a symmetric description of the cone
we write $\delta_{0,\{i\}}$ for $-\psi_i$ in
$\pic {\mgn}.$ (so $\delta_{0.I}$ is
defined whenever $|I|\ge 1$). 

\proclaim{2.1 Theorem} Consider the divisor
$$
D = a \lambda - \birr \dirr - \sum_{
\Sb [g/2] \geq i \geq 0 \\ I \subset N \\ |I| \geq 1 \text{ for }
i = 0 \endSb} b_{i,I} \delta_{i,I}
$$
on $\mgn$ (with the convention that we omit for a given $g,n$ any
terms for which the corresponding boundary divisor does not exist).
Consider the inequalities
\roster
\item $a -12 \birr + b_{1,\emptyset} \geq 0$
\item $\birr \geq 0$
\item $b_{i,I} \geq 0$ for $g-2 \geq i \geq 0$
\item $2\birr \geq b_{i,I}$ for $ g-1 \geq i \geq 1$
\item $b_{i,I} + b_{j,J} \geq b_{i+j,I \cup J}$ for $i,j \geq 0$,
$i + j \leq g-1$, $I \cap J = \emptyset$
\item $b_{i,I} + b_{j,J} + b_{k,K} + b_{l,L} \geq
b_{i+j,I \cup J} + b_{i+k,I \cup K} + b_{i+l, I \cup L}$
for $i,j,k,l \geq 0$, $i + j + k + l =g$, $I,J,K,L$ a partition of $N$
\endroster
where $b_{i,I}$ is defined to be $b_{g-i,I^c}$ for
$i > [g/2]$.

For $g \ge 3$, $D$ has non-negative intersection with all
$1$-dimensional strata iff each of the above inequalities
holds.

For $g =2$, $D$ has non-negative intersection will all
$1$-dimensional strata iff (1) and (3-6) hold.

For $g=1$, $D$ has non-negative intersection with all
$1$-dimensional strata iff (1) and (5-6) hold.

For $g=0$, $D$ has non-negative intersection with all
$1$-dimensional strata iff (6) holds. \endproclaim
\demo{Proof} (2.2) below lists the numerical
possibilities for a stratum. Each inequality above comes
from standard intersection formulae, see e.g.  \cite{Faber96}, 
\cite{Faber97}, 
by intersecting with the corresponding curve of (2.2). Intersecting
with (2.2.4) gives the inequality
$2 \birr \geq b_{i+1,I}$, which shifting indices and notation gives
(2.1.4). \qed \enddemo

Theorem (2.2) gives a listing of numerical possibilities
for one dimensional strata, giving explicit representatives
for each numerical equivalence class. The parts in Theorems (2.1) and
(2.2) correspond: that is, for each family
$X$ listed in parts 2 to 6 of (2.2), the inequality which
expresses the condition that a divisor
$D$ given as in (2.1) meet $X$ non-negatively is given in the
corresponding part of (2.1).

We obtain these
just as in \cite{Faber96}, by defining a map
$\vmgn 0. 4. \rightarrow \mgn$ by attaching a fixed 
pointed curve in some prescribed way. For a subset $I \subset N$,
by a $k + I$ pointed curve we mean a $k + |I|$ pointed curve,
where $|I|$ of the points are labeled by the elements of $I$.

\proclaim{2.2 Theorem} Let $X \subset \mgn$ be a one
dimensional stratum. Then $X$ is either 
\roster
\item For $g \ge 1$, a family of elliptic tails;
\endroster 
or, numerically equivalent to the
image of $\vmgn 0.4. \rightarrow \mgn$ defined by one of the attaching
procedures 2-6 below. 
\roster
\item[2] For $g \geq 3$. Attach a fixed $4 + N$ pointed curve of genus
$g-3$.
\item For $g \geq 2$, $I \subset N$, $g-2 \geq i \geq 0$, $|I| + i > 0$.
Attach a fixed $1 + I$-pointed
curve of genus $i$ and a fixed $3 + I^c$ pointed curve of genus
$g -2 -i$. 
\item For $g \geq 2$ , $I \subset N$, $g-2 \geq i \geq 0$. 
Attach a fixed
$2 + I$-pointed curve of genus $i$ and a fixed $2 + I^c$-pointed
curve of genus $g-2 -i$. 
\item For $g \geq 1$, $I \cap J = \emptyset$, $I,J \subset N$,
$i + j \leq g-1, i,j \geq 0$, $|I| + i,|J| + j > 0$.
Attach a fixed $1 + I$-pointed curve of genus $i$, a fixed $1 + J$-
pointed curve of genus $j$, and a fixed $2 + (I \cup J)^c$-pointed
curve of genus $g-1 -i - j$.
\item For $g \geq 0$,$[I,J,K,L]$ a partition of $N$ into
disjoint subsets, $i,j,k,l \geq 0$,
$i + j + k + l = g$, and 
$i + |I|, j + |J|,k + |K|, l + |L| > 0$. 
Attach $1 + I$,$1+J$,$1+K$, and $1+L$ pointed
curves of genus $i,j,k,l$ respectively.
\endroster

Figure (2.3) shows schematic sketches of each of these 5 families
numbered as in (2.2). 
The generic fiber is shown on the left of each sketch and the 
3 special fibers (with, up to dual graph isomorphism,  any multiplicities) are
shown on the right. The bolder curves are the component(s) of the
$\vmgn 0. 4.$ piece.  Boxes give the type (i.e. genus and marked point set)  
of each fixed component and of each node 
{\it not} of irreducible type. 

In any of the families, if the curve sketched
is not stable, we take the stabilization. Thus e.g. in (6) 
for $n=4$, $g=0$, the map $\vmgn 0.4. \rightarrow \vmgn 0. 4.$
is the identity, while in (3) for $g=2, n =0$, $i = j =0$, the image
of a generic point of $\vmgn 0. 4.$ is the moduli point of 
an irreducible rational 
curve with two ordinary nodes.
\endproclaim

Strata of type (6) play distinctly different roles, both
geometrically and combinatorially, 
from those of type (1-5). Those of type (6) come from
the flag locus, they are the only strata in
genus $0$, and there is a face of $\mc {\mon}.$ that
contains exactly these strata, see (4.9). The cone generated
by the strata (1-5) has a nice geometric meaning given in (0.6).
Our proofs make no direct use of strata of type (6); in particular
we will only directly use the inequalities (2.1.1 -2.1.5).

\draw{110}{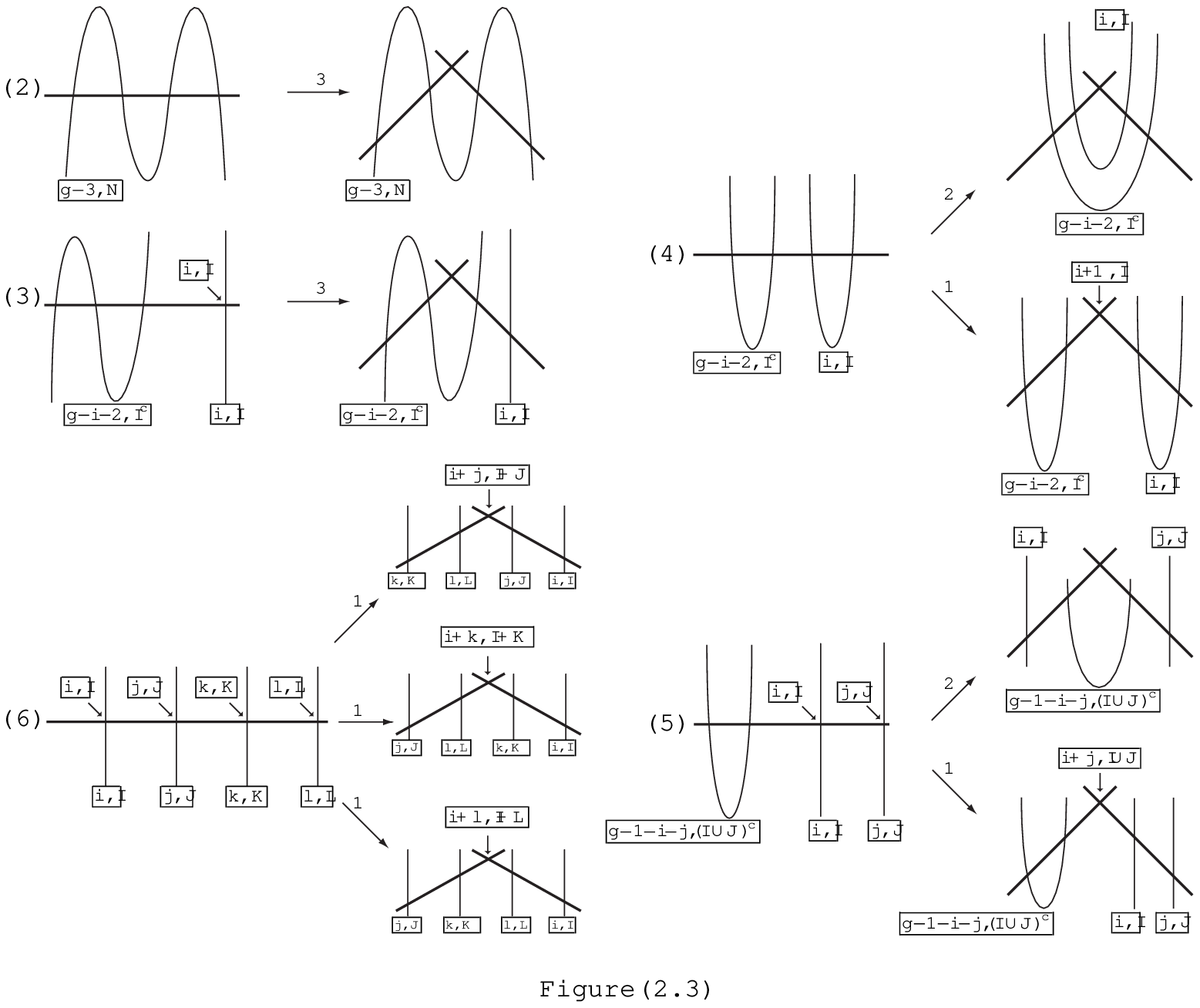}{}%

Strata of type (2-6) are all (numerically equivalent
to) curves lying in $\rgn$.
Those of type (6) are distinguished geometrically 
by the fact that the
curve corresponding to a general point has only disconnecting
nodes, and algebraically by the fact that the corresponding
inequality in (2.1) has more than $3$ (in fact $7$) terms. 

\demo{Proof of 2.2} 
The proof is analogous to that of the case of $n=0$ treated
in \cite{Faber96}. Here are details.

Let $X \subset \mgn$ be a one dimensional stratum: 
$X$ is a 
finite image of a product of moduli spaces (one for each
irreducible component of the curve $C$ corresponding to its general point).  All but one
of these spaces are zero dimensional, thus $\vmgn 0. 3.$, and
there is a distinguished one dimensional factor which is either
$\vmgn 1.1.$ or $\vmgn 0. 4.$. In the $\vmgn 1.1.$ case,
$X$ is a family of elliptic tails and gives (2.1.1) exactly as in \cite{Faber96}. 

Consider
the $\vmgn 0. 4.$ case. Let $C$
be the stable pointed curve corresponding to a general point of $X$ and
$E \subset C$ be the moving irreducible
component (which is rational). Let $h:\tilde{E} \rightarrow E$
be the normalization and let $M \subset \te$ be the union of the
following disjoint subsets:
$h^{-1}(Z \cap E)$ for each 
connected component
$Z$ of the closure of $C \setminus E$, $h^{-1}(p)$ for
each singular point $p$ of $E$, and $h^{-1}(p_i)$ for each
labeled point $p_i$ of $C$ on $E$. These subsets partition
$M$ into equivalence classes. 

We assign to each
equivalence class $\cz$ a triple $(a(\cz),S(\cz),h(\cz))$:
if $\cz = h^{-1}(Z \cap E)$, then 
$h$ is the arithmetic genus of $Z$,
$S \subset \otn$ is the collection of labeled points of $C$ lying on
$Z$, and $a$ is the cardinality of $Z \cap E$. 
$t(h^{-1}(p))=(2,\emptyset,0)$ for a singular point $p$ of $E$, and
$t(h^{-1}(p_i))=(1,\{i\},0)$ for a labeled point $p_i$ of $C$ on $E$. 
The numerical class of $X$ is determined by the collection of triples $\cz$
 by the standard intersection product formulae in, for example, \cite{Faber97}.  
Note that
$\sum_{\cz} ( h(\cz) + a(\cz) -1 ) = g$, that $\bigcup_{\cz} S(\cz)$ is a partition
of $\otn$, and that $\sum_{\cz} a(Z) =4$. 

Now (2-6) are obtained
by enumerating the possibilities
for the collection of $a(Z)$ which are 
$\{4\}$, $\{3,1\}$, $\{2,2\}$, $\{2,1,1\}$ and $\{1,1,1,1\}$.
These correspond to (the obvious generalizations to pointed
curves of) the
families (B), (C), (D), (E), and (F) of \cite{Faber96}
and yield (2-6) by considering the possibilities for
$h(\cz)$ and $S(\cz)$ subject to the above constraints.
\qed \enddemo

Throughout the paper we will separate the $\psi$ and $\delta$
classes (because although they have similar combinatorial
properties, they are very different geometrically, e.g.
$\psi_i$ is nef, while the boundary divisors have very
negative normal bundles). Thus, we will write a divisor as 
$$
\sum_{i=1}^{n} c_i \psi_i + a \lambda - \birr \dirr -
\sum_{\Sb S \subset 
N \\ n-2 \ge |S| \ge 2 \endSb} b_{0,S} \delta_{0,S} - \sum_{[g/2] \geq j
\geq 1,S
\subset
N} b_{i,S} \delta_{i,S} \tag{2.4}
$$
where $b_{0,\{i\}}$ has become $c_i$. As examples, in this notation
 (2.1.5) with $i=j=0,I=\{i\}$ becomes
$$
c_i + b_{0,J} \geq b_{0,J \cup \{i\}}
$$
which in turn specializes to 
$$
c_i + c_j \geq b_{0,\{i,j\}}
$$
in case $J = \{j\}$ while (2.1.3) for $(i,I)=(0,\{i\})$ becomes
$c_i \geq 0$. 

For $g \geq 3$, the expression above is unique (i.e. the various
tautological classes are linearly independent), but for smaller
genera there are relations which are listed in \cite{ArbarelloCornalba98}.
For $g=1$, the boundary divisors are linearly independent
and span the Picard group so 
we can assume above that $a$ and all $c_i$ are zero and, if we do so,
the resulting expression is unique.

\subhead \S 3. Nef on the boundary implies nef \endsubhead

In this section, we prove:

\proclaim{3.1 Proposition} If $g \geq 2$, or $g=1,n\geq 2$, a
 divisor $D \in \tpic \mgn.$  is
nef iff its restriction to $\partial \mgn$ is nef.
\endproclaim

\proclaim{3.2 Lemma} $g \geq 2$. Let $D \in \tpic \mgn.$ be a 
divisor expressed as in (2.1). For $ [g/2] \geq i \geq 1$ 
let 
$$
b_i = \underset {S \subset N} \to {\hbox{\rm max }} b_{i,S}
$$
and define $A \in \tpic \mg.$ by
$$
A = a\lambda - \birr \dirr - \sum_{[g/2] \geq i \geq 1} b_i \delta_i.
$$
If the coefficients of $D$ satisfy one of the inequalities 
(2.1.1-2.1.5) for fixed indices $i,j$ and all subsets of $N$,
then the coefficients of $A$ satisfy the analogous
inequality. In the example of (2.1.5), $b_i + b_j \geq b_{i+j}$ if
$
b_{i,I} + b_{j,J} \geq b_{i+j,I \cup J}$
 for all $i$ and $j$ between $0$ and $ g-1$ and { for all subsets }
$I,J \subset N ${ with } $I \cap J = \emptyset
$.
\endproclaim
\demo{Proof} This is clear since the $b_i$ are defined by
taking a maximum. For example, suppose $b_{i+j} = b_{i+j,T}$.
Then, 
$\displaystyle{
 b_i + b_j \geq b_{i,\emptyset} + b_{j,T} \geq b_{i+j,T} = b_{i+j}} $.\qed
\enddemo

Since all the strata are contained in $\partial \mgn$, (3.1)
will follow immediately from:

\proclaim{3.3 Lemma} Let $D \in \tpic {\mgn}.$. 
If $g \geq 3$ [resp. g=2; g=1]
and $D$ has non-negative intersection with all strata
 in (2.2) of types (1-5) [resp. (1-2) and (4-5); (1) and (5)]
then $D$
is nef outside of $\partial \mgn$, and furthermore
if $g=1$ then $D$ is numerically equivalent to an
effective sum of boundary divisors.
\endproclaim
\demo{Proof} We express $D$ as in (2.4).
We consider first the case of $g=1$. We can as remarked above
assume $a = c_i =0$ for all $i$. Then it is immediate from
(2.1.1) and (2.1.5) (using the translations after (2.4)) 
that all the coefficients are non-positive,
so $D$ is linearly equivalent to an effective sum of boundary
divisors.

Now assume $g \geq 2$. Define $b_i$ and $A$ as in (3.2).

We use the relation
$$
\psi_i = \omega_i + \sum_{i \in S} \delta_{0,S}
$$
(which, as in \cite{ArbarelloCornalba96}, is obtained easily 
from the definition by induction
on $n$)
to write 
$$
 D = \sum_{i \in N} c_i \omega_i +
\biggl(\sum_{\Sb S \subset N \\ n-2 \ge |S| \ge 2\endSb} \Bigl(
\bigl(\sum_{i \in S} c_i\bigr) - b_{0,S}\Bigr)
\delta_{0,S}\biggr) + \pi^*(A) + E \tag{3.4}
$$
where $\pi$ is the forgetful map to $ \mg$ and
$E$ is a combination of boundary divisors $b_{i,S}, i > 0$ ---
an effective one since we take the maximum in forming $b_i$.

By repeated application of (2.1.5) as translated
below (2.4), the coefficients of 
the $\delta_{0,S}$ in (3.4) are non-negative. By (2.1.3), $c_i \geq 0$.
 Since each $\omega_i$ is nef ---
by, e.g. \cite{Keel99} --- it's enough to show $A$ is
nef outside $\partial \mg$. 

Suppose first $g \geq 3$. Then $A$ satisfies each of (2.1.1-2.1.5)
by (3.2). Thus  $a \geq 10\birr \geq 0$ and $a \geq 5b_i$, so the claim
follows from (3.5) below. Now suppose $g=2$. By (3.2) and (2.2) 
$a -10b_{\irr} \geq 2b_{\irr} - b_1 \geq 0$ so the claim
follows from \cite{Faber96,\S 2}. \qed \enddemo

\proclaim{3.5 Lemma} Let $D = a \lambda - \sum b_i \delta_i$
be a divisor on $\mg$, and let $g: Z \rightarrow \mg$ a 
morphism from a projective variety  such that
$g(Z) \not \subset \partial \mg$. Assume
$g \geq 2$. Consider the inequalities:
\roster
\item $a \geq 0$
\item $g \cdot a \geq (8g+4) \birr$
\item $2g \cdot a \geq (8g+4) b_i$, $ [g/2] \geq i \geq 1$
\endroster
If (1-3) hold then $g^*(Z)$ is a limit of effective Weil divisors
and if strict inequality holds in each case and $g$ is 
generically finite then $g^*(Z)$ is big.
\endproclaim
\demo{Proof} The final statement follows from 
the previous statements using $12 \lambda = \kappa + \delta$ and the
ampleness of $\kappa$. 

Let $\Gamma$ be the divisor of \cite{CornalbaHarris88,4.4}
$$
(8g+4) \lambda - g \dirr - 2g \sum_{i>0} \delta_i.
$$

The inequalities imply
$D = s \lambda + c \Gamma + E$ where $c,s \geq 0$, and
$E$ is an effective sum of boundary divisors. The pullback
$g^*(\Gamma)$ is a limit of 
effective Weil divisors by a slight
generalization of \cite{CornalbaHarris88,4.4} (in
characteristic zero we could apply \cite{Moriwaki98}
directly). When $g(Z)$
is contained in the hyperelliptic locus, their analysis
applies without change, as long as the characteristic
is not two. When $g(Z)$ is not contained in the hyperelliptic
locus, one needs the generalization of their (2.9) to
 smooth $T$ of arbitrary dimension (with their
$e_L({\Cal I})$ replaced by the analogous expression supported
at codimension one points of $T$ and the conclusion being
that the expression is a limit of effective Weil divisors, as
in their (1.1)). Their proof extends after 
obvious modifications. \qed \enddemo

\proclaim{3.6 Theorem} Let $D \in \tpic \mgn.$ be a divisor with 
non-negative intersection with all $1$-strata of type
(2.2.1-2.2.5). Let $C \subset \mgn$ be a complete one dimensional
family of curves, whose generic member is a curve without
moving rational components. Then $D \cdot C \geq 0$. If furthermore
$D|_{\rgn}$ is nef, then $D$ is nef.
\endproclaim
\demo{Proof} 
We induce simultaneously on $g,n$. For $g =0$ or 
$g=1,n=1$ there is nothing to prove. So suppose
$g \geq 2$, or $g=1, n \geq 2$. 

Suppose first that $D$ is nef on all $1$-strata of
type (1-5). Let $C$ be as in the statement. We will
show that $D \cdot C \geq 0$. By (3.3) we may assume
$C \subset \partial \mgn$. 
Suppose $C \subset \Delta_{i,S}$ for $i > 0$. We 
apply induction using (1.1), following that notation.
It's enough to show $D_i \cdot C_i \geq 0$. 
A stratum of type (1-5) on one of the factors yields (by the usual
inclusion given by gluing on a fixed curve) a stratum of type
(1-5) on $\mgn$, so $D_i, C_i$ satisfy the induction hypothesis. 
If $C \subset \Delta_{\irr}$ we can pass to $\tilde{\Delta} _{\irr}$
and apply induction directly. 

Now suppose that in addition that $D|_{\rgn}$ is nef. We use
induction just as above (with the same notation) Note that we may assume
that
$D_i$ is nef on $\overline{R}_{ i, S \cup *}$, since we may choose the
fixed curve we glue on to be rational. \qed \enddemo

\subhead \S 4 Reduction to the Flag locus \endsubhead

\proclaim{4.1 Prop} Let $D$ be a divisor on $\mgn$, $g \geq 1$,
satisfying the hypothesis of (3.3).
Let $X \subset \rgn$ be a stratum whose
general member corresponds to a stable curve with
no disconnecting nodes.
Then $D|_{X}$ is linearly equivalent
to an effective sum of boundary divisors and nef divisors.
\endproclaim

\proclaim{4.2 Corollary} 
Let $C \subset \mgn$ be a complete
one dimensional family of pointed curves whose general member
has no moving smooth rational components. 
Let $D$ be a divisor on $\mgn$. If $D$ has non-negative
intersection with all $1$-strata of type (2.2.1-2.2.5), then
$D \cdot C \geq 0$. 
\endproclaim

\demo{Proof} We use induction as in the proof
of (3.6). By (1.2) we can assume that the general member
of $C$ is irreducible, and by (3.6), rational. Now
(4.1) applies. \qed \enddemo

\proclaim{4.3 Corollary} $$ \tomc {\mgn}. = N + \mc {\vmgn 0. g+n.}/S_g. $$
where  $N$ is the subcone generated by the strata (2.2.1-2.2.5).
\endproclaim
\demo{Proof} Immediate from (4.2), (1.1) and (1.2). \qed \enddemo

Of course (4.2) and (4.3) together contain (0.6), which in turn implies (0.3).

\demo{Proof of (4.1)}
First we reduce to the case of $n=0$. As in (3.4) we have
$$
D = \sum_{i \in N} c_i \omega_i
+ \pi^*(A) + G
$$
with $G$ an effective combination
of boundary divisors parameterizing degenerations
with a disconnecting node. (In the genus $1$ case the
sum is empty and $A$ is
a multiple of $\delta_{\irr}$). So we can replace $X$ by its
image under $\pi$. When $g=1$ this is a point and there
is nothing to prove. When $g=2$ this image is either a point,
or is one dimensional and defines inequality (2.1.4) (so the
pullback of $D$ is nef). So we can assume  $n=0, g \geq 3$.

Let $C$ be the stable curve corresponding to a general
point of $X$, and $E$ a component of $C$.
As in the proof of (2.2) we get from $E$ a map of some $\mon$ 
into $\mg$ and we
want to see the pullback of $D$ is an effective
sum of boundary components and nef divisors.
As in the proof of (2.2), the $n$ points are
divided into equivalence classes: the singular points
of $E$ are divided into disjoint $2$-cycles, while
the points where $E$ meets
other irreducible components of $C$ are divided by the connected components
of (the closure of) $C \setminus E$, each class with at
least $2$ elements by our assumption that $C$ has no disconnecting nodes.
If $P$ is the corresponding partition, then a class $\delta_i$, $i > 0$
pulls back to a sum of irreducible divisors $\delta_{\ct}$ with $\ct > P$
(i.e. $P$ refining the partition $[T, T^c]$). We use the identity
on $\mg$
$$
- \dirr + 12 \lambda = \kappa + \sum_{i>0} \delta_i
$$
and the fact that $\lambda$ is trivial on $\mon$. If 
$b_{\irr} \leq 0$ then $b_{\irr} =0$ by (2.1.2) and
so the result is clear from (2.1.4).
So we may assume (by scaling) that $b_{\irr} =1$. Then 
$D$ pulls back to
$$
\kappa + \sum_{\ct > P} e_T \delta_{\ct}
$$
with $e_T \geq -1$ by (2.1.4). Now apply
(4.4) below.
\qed
\enddemo

\proclaim{4.4 Lemma} Let $P$ be a partition of $\otn$, such that
each equivalence class contains at least two elements.
The divisor
$$
\kappa + \sum_{\ct > P} e_{\ct} \delta_{\ct} \tag{4.5}
$$
on $\mon$ 
is linearly equivalent to an effective sum of boundary divisors,
if each $e_{\ct} \geq -1$.
\endproclaim
\remark{Remark} Note that at least some condition on $P$ is
required. E.g. $\kappa - \delta$ is not an effective sum
of boundary divisors, since as a symmetric expression the
coefficient of $B_2$ 
is negative, and there are no symmetric relations.
\endremark
\demo{Proof of (4.4)}
Expanding $\kappa$ in terms of boundary divisors (see e.g.
\cite{KeelMcKernan96}) we have
$$
\sum_{\ct > P} (\frac{|T|(n-|T|)}{n-1} + e_{\ct} -1) \,\delta_{\ct} +
\sum_{ \ct \not > P} (\frac{|T|(n-|T|)}{n-1} -1) \, \delta_{\ct}.
$$
Note $\frac{|T|(n-|T|)}{n-1} \geq 2$ unless
$|T| =2$, in which case it is $2 - \frac{2}{n-1}$, or $n=6$
and $|T| = 3$ in which case it is $2 - \frac{1}{n-1}$. So the sum is
nearly effective as it stands. We use relations in the Picard
group to deal with the negative coefficients.

We assume first $n \geq 7$.
If there are no $2$-cycles in $P$ our sum is already effective, so 
 (after perhaps reordering) we can write $P$ as
$$
(12)(34)\dots(2k-1 \  2k)
$$
for some $k \ge 1$.
Further we can assume either that $2k=n$, or we have
a further $t$-cycle, $t \geq 3$, $(2k+1 \  2k+2 \ \dots\ 2k+t)$.
We use the $k$ 
relations
$$
\sum_{[\{i,i+1\},\{i+2,i+3\}] = \ct|_{\{i,i+1,i+2,i+3\}}} \delta_{\ct}
= \sum_{[\{i,i+2\},\{i+1,i+3\}] = \ct|_{\{i,i+1,i+2,i+3\}}} \delta_{\ct} 
$$
for $i = 1,3,5, \dots, 2k-1$
with the adjustment that if $2k =n$, then we replace the last
by the analogous relation for $[\{2k-1,2k\},\{1,2\}]$. By
$[\{i,i+1\},\{i+2,i+3\}] = \ct|_{\{i,i+1,i+2,i+3\}}$ we mean
that the restriction of $\ct$ to the $4$-element set is the
given partition into disjoint $2$-cycles, i.e. 
either $i,i+1 \in T$ and $i+2, i+3 \in T^c$ or this holds
with $T$ and $T^c$ reversed. Thus the above relations come directly
from those in \cite{Keel92} by adjusting notation.

Let $L = R$ be the relation among effective boundary divisors that one
obtains by adding together the 
$k$ relations above (without canceling terms common to both sides).
Let
$\dirr$ be the sum of
the $\delta_{\ct}$ with $\ct \not > P$, and $\delta_1$ the sum
of those with $\ct$ one of the $2$-cycles of $P$
together with its complement. Note
$R$ is supported on $\dirr$ while $L \geq \delta_1$ (here
and throughout the paper inequalities between sums of Weil divisors stand
for a system of inequalities, one for each component) and $L \geq 2
\delta_1$ when $n=2k$.

Choose $s>0$ maximal so
that
$$
s R \leq \sum_{\ct \not > P} (\frac{|T|(n-|T|)}{n-1} -1)\,  \delta_{\ct}.
$$
Thus (4.5) is linearly equivalent to
$$
sL  - \frac{2}{n-1} \delta_1 + E
$$
for $E$ an effective sum of boundary divisors.

Define coefficients $\alpha_{\ct}$ so
$R = \sum_{\ct} \alpha_{\ct} \delta_{\ct}$.

Since $L \geq \delta_1$, it is enough to show that
$s\geq \frac{2}{n-1}$, or equivalently that
$$
\frac{2}{n-1} \alpha_{\ct}  \leq \frac{|T|(n-|T|)}{n-1} -1 \tag{4.6}
$$
and for $n =2k$, $\frac{\alpha_{\ct}}{n-1} \leq \frac{|T|(n-|T|)}{n-1} -1$
is sufficient.
Suppose first $|T| =2$. Note that $\delta_{0,T}$ can occur
(with non-zero coefficient) on the right hand side of
at most one of the $k$ relations. Thus $\alpha_{\ct} \leq 1$ and
(4.6) follows.

Suppose now
$|T| \geq 3$. Then (4.6) holds as long as $2 \alpha_{\ct} \leq n-1$.
Since $\alpha_{\ct} \leq k$, we can assume $2k =n$. Then we only
need $\alpha_{\ct} \leq 2k-1$, which holds as $\alpha_{\ct} \leq k$.

Now we consider the cases of $n \leq 6$. If $P$ has an $n$-cycle,
the result is obvious ($\ct > P$ is then impossible). Up to reordering the
possibilities are  $(12)(34)$, $(12)(345)$,
$(12)(34)(56)$ and $(123)(456)$. In the first case (4.5)
is nef (and we are working on $\vmgn 0.4. = \pr 1.$). In the second
case the above argument applies without change. In the third we
can combine the above argument for the cases $n =2k$ and $|T|=2$.
Finally consider $P = (123)(456)$.
Thus
$\ct > P$ implies $\delta_{\ct} = \delta_{0,{\{1,2,3,\}}}$
We use the relation
$$
\sum_{[\{1,2\},\{5,6\}] = \ct|_{\{1,2,5,6\}} } \delta_{\ct} =
\sum_{[\{1,5\},\{2,6\}] = \ct|_{\{1,2,5,6\}} } \delta_{\ct}
$$
to show (4.5) is effective as above. \qed \enddemo

\proclaim{4.7 Theorem} Consider the line bundle
$$
D = a \lambda + \sum_{i=1}^{i=n} (g+n-1)\psi_i -
\birr \dirr - \sum \bigl(g+n-(i + |S|)\bigr)\bigl(i + |S|\bigr) \delta_{i,S}
$$
on $\mgn$.

Then, $D|_{\fgn}$ is trivial and $D$ has zero intersection
with all strata of type (2.2.6). Moreover, $D$ has
strictly positive intersection with all other strata 
if
$a > 12\birr -(g+n-1)$, $2\birr > ((n+g)/2)^2$. If $D$
satisfies these inequalities, $D$ is nef.
In particular $F_1(\mgn)$ implies $F_1(\vmgn 0. g+n. /S_g)$.

Every nef line bundle on $\vmgn 0. g+n./S_g$ is the pullback
of a nef line bundle on $\mgn$. The strata numerically equivalent
to a curve in the flag locus are precisely the strata of type (2.2.6).
\endproclaim
\demo{Proof} One checks by (2.1.6) that $D$ has trivial intersection
with all the strata of type (2.2.6). These correspond to the
strata of $\vmgn 0. g + n.$. 
Since
the strata in $\mon$ generate the Chow group, it follows
that $D$ has trivial pullback, as in (0.3), to $\vmgn 0. g+n. /S_g$.  
Now one checks using (2.1-2.2) that
$D$ has strictly
positive intersection with all other strata if
the inequalities in the statement hold. The final remark in
the statement now follows.

Now suppose the inequalities of the statement hold. Then
$D$ is nef by (0.3). Standard intersection formulae,
e.g. \cite{Faber97}, show that 
$$
f^*: \tpic \mgn. \rightarrow \pic {\vmgn 0. g+n./S_g}.
$$ is surjective.
Suppose $f^*(W)$ is nef. Then $W + mD$ will have
non-negative intersection with all strata for $m > > 0$
by the above. Thus it is nef by (0.3). \qed \enddemo

Of course the above theorem implies (0.7).

\proclaim{4.8 Lemma} Let $D= D_{g,n}$ be the divisor of (4.7). Let
$h: \vmgn j. s. \rightarrow \mgn$ be a factor of one of 
the strata. If $h$ corresponds to a smooth rational component
then $h^*(D)$ is trivial. 
Otherwise $\ex h^*(D). \subset \partial \vmgn j. s.$. \endproclaim
\demo{Proof} In the first case we can assume (without
affecting the pullback on Picard groups) that $h$ comes from
a stratum in the flag locus, and so the result follows from
(4.7). 

Now assume either $j \geq 1$, or $j=0$ and the corresponding
component contains at least one disconnecting node. If 
$s \geq 1$, let
$C$ be the general fiber of one of the projections
$\vmgn j. s. \rightarrow \vmgn j. s-1.$. Then, $C$ is numerically
equivalent to an effective sum of strata, not all of which
are of type (2.2.6). Thus $h^*(D) \cdot C > 0$ by (4.7). Now
(0.10) applies as long as $j \geq 2$.

Suppose $j=1$. Observe by (4.7) that $d D_{j,s} - h^*(D)$
is nef for sufficiently large $d$. Thus its enough to show
$\ex D_{j,s}. \subset \partial \vmgn j. s.$. By (4.7) and
the proof
of (3.3), $D_{1,s}$ is an effective sum of boundary divisors,
with every divisor occurring with positive coefficient. Since
the boundary supports an ample divisor, the result follows.

Now assume $j=0$, and the corresponding component contains
a disconnecting node. As in the previous paragraph, we
can assume $h$ defines $\rgn$, $g \geq 1$. By (4.7),
for any divisor $W$, $dD - W$ satisfies the hypothesis of (3.3).
So, taking $W$ ample, the result follows from (4.1). \qed \enddemo

\proclaim{4.9 Corollary} If $D$ is the nef divisor of (4.7), then
$$
D^{\perp} \cap \mc {\mgn}. = \mc {\vmgn 0. g+n.}/S_g. .
$$
\endproclaim
\demo{Proof} Let $C \subset \mgn$ be a curve. By (4.8) and 
(1.2), $D \cdot C = 0$ iff every moving component of the 
associated family of curves is smooth and rational, in which
case $C$ is algebraically equivalent to a curve in $\fgn$. 
\qed \enddemo

\subhead \S 5 Exceptional Loci  \endsubhead

In this section we prove (0.9). We will use the following lemmas.

\proclaim{5.1 Lemma} Let $F$ the fiber of
$\pi: \mgn \rightarrow \mg$ (resp.
$\pi: \vmgn 1.n. \rightarrow \vmgn 1.1.$) over
a point $[C]$ (resp. $[C,p]$) with $C$ smooth. Let 
$G$ be the automorphism group of $C$. Let $D$ be
divisor on $F$ which is the restriction of a nef
divisor in $\tpic \mgn.$. Let 
$$
f: F \rightarrow C^n/G
$$
be the natural birational map. $D$ is big iff
$f_*(D)$ is big, in which case
$$
\ex D. \subset f^{-1}(\ex {f_*(D)}.) \cup \partial F.
$$
\endproclaim
\demo{Proof} We will abuse notation slightly and use the same
symbol for a divisor on $\mgn$ and its restriction to~$F$.

Suppose first $g \geq 2$. As in (3.4) we have
an expression
$$
 D = \sum_{i \in N} c_i \omega_i +
\biggl(\sum_{\Sb S \subset N \\ n-2 \ge |S| \ge 2\endSb} \Bigl(
\bigl(\sum_{i \in S} c_i\bigr) - b_{0,S}\Bigr)
\delta_{0,S}\biggr) 
$$
with $c_i \geq 0$, and with each $\delta_{0,S}$ having non-negative
coefficient. Let
$$
\Gamma =\biggl(\sum_{\Sb S \subset N \\ n-2 \ge |S| \ge 2\endSb} \Bigl(
\bigl(\sum_{i \in S} c_i\bigr) - b_{0,S}\Bigr)
\delta_{0,S}\biggr).
$$
Observe that the exceptional divisors
of $f$ are precisely the $\delta_{0,S}$ with
$|S| \geq 3$, and that 
$f(\delta_{0,S}) \subset f(\delta_{0,T})$
iff $T \subset S$. Thus the support of
$f^{-1}(f_*(\Gamma))$ is the union of the
$\delta_{0,S}$ over subsets $S$ which contain
a $2$-element subset, $T$, such that $\delta_{0,T}$
occurs in $\Gamma$ with non-zero coefficient.
By repeated application of (2.1.5), we obtain
$$
(\sum_{i \in S} c_i) - b_{0,S} \geq 
(\sum_{i \in T} c_i ) - b_{0,T}
$$
if $T \subset S$. Thus the support of $\Gamma$ contains
the support of 
$f^{-1}(f_*(\Gamma))$ (in fact since $D$ 
is $f$-nef, by negativity
of contractions, the supports are the same). 
As each $\omega_i$
is a pullback along $f$ it follows 
that $D$ is big iff $f_*(D)$ is big.
By \cite{Keel99,4.9}, $\omega_i$
is nef and $\ex \omega_i. \subset \partial F.$ 
The result follows.

Now suppose $g=1$. 
By (3.3) we may express $D$ as 
$$
-\sum b_{0,S} \delta_{0,S} 
$$
with $b_{0,S} \leq 0$. Using (2.1) one checks
$b_{0,S} \leq b_{0,T}$ for $T \subset S$.
Now it follows as above that $D$ has the same support as
$f^{-1}(f_*(D))$. The result follows. \qed \enddemo

\proclaim{5.2 Lemma} Let $C$ be a smooth curve of genus
$g \geq 2$, with automorphism group $G$. Let
$D$ be an divisor on $C^n/G$ of the form
$$
D = \sum c_i \omega_i + \sum_{i \neq j} a_{ij} \Delta_{ij}
$$
with $c_i, a_{ij} \geq 0$ and $c_i + c_j \geq a_{ij}$.
$D|_Z$
is big for any subvariety $Z \subset C^n/G$ not contained
in any diagonal iff 
$D$ is $q_i$-relatively ample for each of the  projections
$$
q_i: C^n/G \rightarrow C^{n-1}/G
$$
given by dropping the $i^{th}$ point.
\endproclaim
\demo{Proof} The forward implication is obvious.
Consider the reverse implication. 
As $G$ is finite it's enough to consider the
analogous statements on $C^n$. We assume $D$ is 
$q_i$-ample for all $i$. Write
$$
D = \sum c_i \omega_i + \sum a_{ij} \Delta_{ij}.
$$
Reorder so that $c_i = 0$ iff $i \geq {s+1}$.
For each $t>s$, since $D$ is $\pi_t$-ample,
$a_{i(t) t} > 0$  for some $i(t) \neq t$. 
Since $c_i + c_j \geq a_{ij} \geq 0$, we must have $i(t) \leq s$. 
If
$$
p: C^n \rightarrow C^{s}
$$ 
is the projection onto the
first $s$ factors, then 
$\sum c_i \omega_i$ is the pullback along $p$ of an
ample divisor, while 
$$
\sum_{t > s} a_{i(t) t} \Delta_{i(t) t}
$$
is $p$-ample. \qed \enddemo

\proclaim{5.3 Lemma} Let $C$ be a smooth curve of genus $1$,
with automorphism group $G$. Let
$$
D = \sum a_{ij} \Delta_{ij}
$$
be an effective combination of diagonals on $C^n/G$. 
Let $H \subset S_n$ be the subgroup generated by
the $(ij)$ with $a_{ij} \neq 0$. If $H$ is transitive,
then $D$ is ample. If a partition $N = S \cup S^c$
is preserved by $H$, then $D$ is pulled back along the
fibration
$$
C^n/G \rightarrow C^S/G \times C^{S^c}/G.
$$
\endproclaim
\demo{Proof}   
The second claim is clear. For the first we induce
on $n$. For $n=2$ the result is clear
(the quotient is a curve and $D$ is non-trivial).
Note since the genus is one, that
all the diagonals (and hence $D$) are semi-ample, so its enough to
show $D$ has positive intersection with every irreducible
curve. 
But if $\Delta_{ij} \cdot E = 0$ for a curve $E$,
then $E$ is numerically equivalent to a curve contained
in $\Delta_{ij}$ so it's enough to show
$D|_{\Delta_{ij}}$ is ample, and for this
we can apply induction. \qed \enddemo

\demo{Proof of (0.9)} 
Suppose first $g \geq 2$. 
If  $n=0$ and $g=2$, the
result follows from \cite{Faber96}.
If $n=0$ and  $g \geq 3$, we
apply (3.5). If $a =0$, then $D$ is trivial by (2.1.1-2.1.4).
Otherwise we have strict inequality in (3.5) by (2.1).

So we may assume $n > 0$ and that $D$ is not pulled back
from $\vmgn g. n-1.$. If $D$ is trivial
on the general fiber of $q: \mgn \rightarrow \vmgn g. n-1.$, 
then it is $q$-trivial, since it is nef, 
thus, by an
easy calculation, pulled back. So we can assume $D$ 
is $q$-big (for any choice of $q$).

Mimicking the argument of (3.3) we write $D$ as in (3.4)
$$
D = \bigl(\sum_{i \in s} c_i \omega_i \bigr) +
\biggl(\sum_{\Sb S \subset \otn \\ n-2 \ge |S| \ge 2\endSb} \Bigl(
\bigl(\sum_{i
\in S} c_i\bigr) - b_{0,S}\Bigr)
\delta_{0,S}\biggr) + \pi^*(A) + E.
$$
For any subvariety $Z \not \subset \partial
\mg$, $A|_Z$ is effective by (3.5). 
Thus by \cite{Koll\'arMori98.3.23} it is enough to show
that $\ex {D|_F}. \subset F \cap \partial \mgn$ for
$F$ as in (5.1). This is immediate from (5.1-5.2).

Now suppose $g=1$. Suppose first that $D$ is
$q_S$ big for each fibration
$$
q_S: \vmgn 1.n. \rightarrow \vmgn 1. S. \times_{\vmgn 1.1.} \vmgn 1. S^c. 
$$
and each non-trivial proper subset $S \subset N$.
We conclude as in the $g \geq 2$ case by applying (5.1)
and (5.3).

So we may assume we have a non-trivial proper subset $S \subset N$
such that $D$ has zero intersection with the general fiber (which is
one dimensional) of $q_S$. 
As $D$ is nef, a short calculation shows $D$ is of the
form $\pi_S^*(D_S) + \pi_{S^c}^*(D_{S^c})$ for
divisors $D_S$ and $D_{S^c}$ on $\vmgn 1.S.$ and $\vmgn 1.S^c.$.
Indeed 
one checks with
(2.2.5) that 
$$
b_{0,T \cap S} + b_{0,T \cap S^c} = b_{0,T}
$$
for all $T$, from which it follows that we
may take 
$$
\align
D_S &= -\sum_{T \subset N} b_{0,T \cap S} \delta_{0,T \cap S} \\
D_{S^c} &= -\sum_{T \subset N} b_{0,T \cap S^c} \delta_{0,T \cap S^c}.
\endalign
$$
Now note that $\vrgn 1. m.$ and all of the $1$-strata
of $\vmgn 1.n.$ except for the elliptic tails lie over
the point $\dirr \in \vmgn 1.1.$. It follows from
(3.6) that for divisors $D_S \in \pic {\vmgn 1.S.}.$ and
$D_{S^c} \in \pic {\vmgn 1. S^c.}.$, if 
$\pi_{S}^* D_{S} + \pi_{S^c}^*D_{S^c}$ is nef, then both
$D_{S}$ and $D_{S^c}$ are nef, as long as each has non-negative
intersection with a family of elliptic tails. 
Let $E \subset \vmgn 1. n.$ be a family of elliptic curves.
Then so are  $E_S=\pi_S(E)$ and $E_{S^c} =\pi_{S^c}(E)$ and 
$\dirr \cdot E_S = \dirr \cdot E_{S^c}> 0$.
So if we define $t$ by
$$
(D_S + t  \dirr) \cdot E_S =0
$$
then 
$$
(D_{S^c} - t\dirr) \cdot E_{S^c} 
= D_S \cdot E_S + D_{S^c} \cdot E_{S^c} = D \cdot E \geq 0.
$$
Hence 
$$
D = \pi_{S}^*(D_S + t \dirr) + \pi_{S^c}^*(D_{S^c} - t \dirr)
$$
with both
$D_S + t \dirr$ and $D_{S^c} - t \dirr$
nef. \qed \enddemo

\remark{Remarks} A natural analog of $f$ from
(5.1) for $g=0$ is the birational map 
$$
f: \mon \rightarrow Q
$$
with $Q$ one of the maximal GIT quotients of $(\pr 1.)^{\times n}$
(e.g. the unique $S_n$ symmetric quotient).
The birational geometry of $Q$ is well understood,
see \cite{HuKeel00}, and in particular the analog of
(5.3) holds.  Thus if one could prove the analog
of (5.1), one could extend (0.9) to all genera. This 
would imply (0.2). \endremark

\subhead \S 6 Ad hoc examples \endsubhead

\proclaim{6.1 Proposition} {\rm(char $0$)}
Let $D$ be an  $F$-nef divisor $a \lambda - \sum b_i
\delta_i$ on $\mg$. Assume further
that for each coefficient $b_i$,
$[g/2] \geq i \geq 1$ that either $b_i =0$ or $b_i \geq \birr$.
Then $D$ is nef.
\endproclaim
\demo{Proof} By (0.3) it is enough to show $h^*(D)$ is
nef, for the composition $h: \vmgn 0. 2g. \rightarrow \rgn \subset
\mg$. We can assume $\birr > 0$ (otherwise all $b_i =0$ ).
Suppose first $b_i = 0$ for some $i > 0$. Consider the
map
$$
\vmgn i. 1. \rightarrow \mg
$$
(from the usual product decomposition of $\delta_i$). Let
$D'$ be the pullback of $D$, which we express in the
usual basis.
Intersecting with an appropriate
stratum of type (2.1.5)
on $\vmgn i. 1.$ gives the inequality
$$
c + b_{j,\emptyset} \geq b_{j,\{1\}} = b_{g-j,\emptyset},\quad g > j > 0
$$
(where $c=c_1$). In our case $c = b_i= 0$ and so $b_{j, \emptyset} = b_{j,
\{1\}}$, and thus $D'$ is pulled back from $\vmgn i. 0.$. The same
argument applies to the other factor of $\delta_i$. Thus by
induction on $g$ we may assume $D|_{\delta_i}$ is
nef for any $i > 0$ with $b_i = 0$.

Using $12 \lambda = \kappa + \delta$ and the expression
$K_{\mon} = \kappa + \delta$ we have that
$$
f^*(D/\birr) = K + \sum_{i =1}^{[g/2]} (2 - b_i/\birr) \delta_i + B
$$
where $\delta_i$ on $\vmgn 0. 2g.$ indicates the pullback
of $\delta_i$ on $\mg$ (note the pullback of $\delta_i$
to $\rgn$ is effective and irreducible) which is an effective sum
of boundary divisors with coefficient one (or zero), and
$B$ is the sum of the boundary divisors of $\vmgn 0. 2g.$
other than the $\delta_i$.
Let $\Delta =  \sum_{i =1}^{[g/2]} (2 - b_i/\birr) \delta_i+ B$.
By (2.1.4) the coefficients of $\Delta$ are non-negative
and, by the above, $K + \Delta$ is nef on any boundary component whose
coefficient in $\Delta$ is greater than one.
Suppose $D$ (or equivalently $K + \Delta$) 
is not nef. Then there is a $K + \Delta$ negative
extremal ray
$R$ of $\mc {\vmgn 0. 2g.}.$. If $\Delta'$ is obtained from
$\Delta$ by shrinking to one the coefficients that are 
greater than one, 
then by the above $(K + \Delta') \cdot R < 0$. Thus by
\cite{KeelMcKernan96,1.2}, $R$ is generated by a one dimensional
stratum, a contradiction. \qed \enddemo

Note (6.1) applies to all the edges of the Faber cone computed in
\cite{Faber96} for $g \leq 7$,
giving another proof of (0.2) for $n=0$, $g \leq 7$.
Unfortunately, this fails for $g=10$ where
$30\lambda -3 \delta_0 - 6 \delta_1 - 6\delta_2 -2 \delta_3 - 4\delta_4
-6\delta_5$ is a vertex. It also answers in the
affirmative question (a) at the end of \S 2 of \cite{Faber96}.

\proclaim{6.2 Corollary} {\rm(char $0$)} The ray $10 \lambda - 2 \delta + \dirr$
is nef on $\mg$ for all $g \geq 2$.
\endproclaim

\noindent Also, by averaging with the nef ray $12 \lambda - \dirr $, it
recovers the ampleness result of  \cite{Cornalba-Harris88}.

\proclaim{6.3 Corollary (Cornalba-Harris)} {\rm(char $0$)} The class $11 \lambda
- \delta$ is nef on $\mg$ for all $g \geq 2$.
\endproclaim

The ray $11 \lambda - \delta$ is known to be semi-ample and the
associated map to be a birational contraction of relative
Picard number one with exceptional
locus $\Delta_1$ which blows down families of elliptic tails. In fact, except
for some special low genus examples, this
is the only elementary divisorial blowdown of $\mg$:

\proclaim{6.4 Proposition} {\rm (char $0$, $g \geq 5$)}
The only divisorial contraction $f: \mgn \rightarrow X$ of relative Picard number
one with $X$ projective is the blowdown of the elliptic tails.
\endproclaim
\demo{Proof} By (0.9) some boundary divisor is $f$-exceptional.
We assume $\delta_i$ is exceptional, for some $i \geq 1$, the
argument for $\delta_{\irr}$ is only notationally different.
By the product decomposition, $f$ induces a fibration on 
$\vmgn k. 1.$  either for $k=i$ or $k = g-i$. Assume $f$ 
does not contract the elliptic tails. Then by (0.9), this 
induced fibration
factors
through the map $\pi_k: \vmgn k. 1. \rightarrow \vmg k.$.
For $i > 1$, the general fibers of $\pi_i$ and $\pi_{g-i}$
generate the same ray in $N_1(\mg)$ (indeed, the only basic 
class with non-zero intersection with either is $\delta_{i}$, which
intersects either negatively). Thus for $i > 1$, $f$ induces a
fibration on either factor, and contracts any curve contracted
by $\pi_k$, $k = i, g-i$. But for $k > 2$ (which holds for 
at least one of $i$ or $g-i$ since $g \geq 5$) the image in 
$N_1(\mg)$ of $N_1(\pi_k)$ (the 
relative Neron-Severi group of curves of $\pi_k$ ) has dimension
at least two. Indeed the fiber of $\pi_k$ over
a general point of $\delta_1(\vmg k.)$ is a curve with
two irreducible components, and a quick computation shows
the images in $N_1(\mg)$ of these components are linearly independent 
($\delta_1$ intersects one negatively, and the other positively).
\qed \enddemo

\remark{Remark} Bill Rulla, \cite{Rulla00} 
has shown that on $\vmg 3.$ there
is a second relative Picard number one contraction 
of $\delta_1$. The corresponding extremal
ray is generated by the curve (2.2.4), with $i=1$.
\endremark

\subhead \S 7 Reformulation of $F_1(\mon)$ and review of the
evidence \endsubhead

Let $S \subset \mon$ be the
union, with reduced structure, of the two dimensional strata. The
irreducible  components of $S$ are all smooth Del Pezzo surfaces,
either $\vmgn 0. 5.$ or $\pr 1. \times \pr 1.$.

\proclaim{7.1 Proposition} The natural map
$A_1(S) \rightarrow A_1(\mon)$ is an isomorphism, and
under this map $\mc S.$ is identified with the Faber cone.
Thus $F_1(\mon)$ is equivalent to the statement that the
natural injection $\mc S. \rightarrow \mc {\mon}.$
is an isomorphism.
\endproclaim
\demo{Proof} The strata generate the Mori
cone of $S$, so it's enough to prove the statement about
Chow groups. This follows from the Chow group description 
given in \cite{KontsevichManin96}.
\qed \enddemo

It is known that each stratum of $\mon$ generates an
extremal ray. Each can be contracted by a birational
morphism of relative Picard number one which is in fact a log Mori
fibration( cf. \cite{KeelMckernan96}). The main evidence
for $F_1(\mon)$ is \cite{KeelMckernan96,1.2}, which states in
particular 
that (0.2) holds for any extremal ray of $\mc {\mon}.$ that can be 
contracted by a morphism of relative Picard number one,
except possibly when the
map is birational and has one dimensional exceptional
locus. Any stratum
deforms once $n \geq 9$, so (0.2) would imply this extra
condition always holds for such $n$. In practice the 
main way of showing a curve generates an extremal way
is to exhibit a corresponding contraction, so if (0.2)
fails, it will be a challenge to find a counter-example.

Kapranov proved that $\mon$ is the inverse limit of
all GIT quotients $Q$ (varying the polarization)
for
$\sltwo$ acting diagonally on $(\pr 1.)^n$. In
particular each $Q$ is a contraction of $\mon$, and 
so inherits a stratification. In 
\cite{HuKeel00} $F_1(Q)$ and $F^1(Q)$ are proved.
One can also consider
the contractions $\mon \rightarrow \vln {n-2}.$, 
introduced in \cite{LosevManin00}. $\vln {n-2}.$ is a
smooth projective toric variety. On any projective
toric variety the toric strata generate the cones
of effective cycles. The toric strata for $\vln {n-2}.$ are images
of strata in $\mon$, thus $F_k(\vln {n-2}.)$ holds for
all $k$. We remark that $\vln {n-2}.$ gives particularly good 
evidence for $F_1(\mon)$
as it has Picard
number $2^{n-2} - n + 1$, roughly half the Picard
number of $\mon$  whereas a GIT quotient $Q$ can
have Picard number at most $n$.

%
%
%
%

%
%
%
%

\allowbreak


\Refs
\widestnumber\key{ArbarelloCornalba96}
\ref \key{ArbarelloCornalba96} \by E.~Arbarello and M.~Cornalba
      \paper Combinatorial and
algebro-geometric cohomology classes on the moduli spaces of curves
       \jour J. Algebraic Geom. 
       \vol 5
       \yr 1996
       \pages 705--749
\endref
\ref \key{ArbarelloCornalba98} \bysame
     \paper Calculating cohomology groups of moduli spaces
of curves via algebraic geometry
     \jour Inst. Hautes \'Etudes Sci. Publ. Math. 
     \vol  88
      \yr 1998
      \pages 97--127
\endref
\ref \key{CHM97} \by L.~Caporaso, J.~Harris, and B.~Mazur
    \paper Uniformity of rational points
    \jour  J. Amer. Math. Soc. 
     \vol 10
     \yr 1997
    \pages 1--35
\endref
\ref \key{CornalbaHarris88} \by  M.~Cornalba and J.~Harris
     \paper Divisor classes associated to families of stable varieties, with
applications to the moduli space of curves
     \jour Ann. Sci. \'Ecole Norm. Sup. (4)
     \vol 21
     \pages 455--475
     \yr 1988
\endref

\ref \key{DeligneMumford69} \by P.~Deligne and D.~Mumford
     \paper The irreduciblility of the space of curves of given genus
      \jour Inst. Hautes \'Etudes Sci. Publ. Math. 
      \vol 36
      \yr  1969 
      \pages 75--109
\endref
\ref \key{Faber96} \by  C.~Faber
      \paper Intersection-theoretical computations on 
$\mg$
      \jour  Banach Center Publ.
      \vol  36, 
      \yr 1996
      \pages 71--81
\endref
\ref \key{Faber97} \bysame
      \paper Algorithms for computing intersection numbers on
moduli spaces of curves, with application to the class of the
locus of Jacobians
      \yr 1999
      \inbook London Math. Soc. Lecture Note Ser.
      \vol 264
      \pages 93--109
\endref
\ref \key{Gibney00} \by  A.~Gibney
     \book Fibrations of $\mgn$
      \bookinfo Ph. D. Thesis, Univ. of Texas at Austin
      \yr 2000
\endref
\ref\key{HarrisMorrison98} \by  J.~Harris and I.~Morrison
    \book Moduli of Curves 
\bookinfo Grad. Texts in Math.
    \publ Springer-Verlag
\publaddr New York
    \vol 187
    \yr 1998
\endref
\ref \key{HuKeel00} \by  Y.~Hu and S.~Keel
     \paper Mori dream spaces and GIT
      \jour Michigna Math. J. 
      \vol 48 
      \yr 2000
      \pages 331--348
\endref
\ref \key{Keel92} \by S.~Keel
     \paper Intersection theory of moduli space of stable $N$-pointed curves of
genus zero 
      \jour Trans. Amer. Math. Soc. 
      \vol 330 
      \yr 1992
      \pages 545--574
\endref
\ref \key{Keel99} \bysame
     \paper Basepoint freeness for nef and big linebundles in positive
characteristic
     \jour Annals of Math.
     \yr 1999
      \pages 253--286
\endref
\ref\key{KeelMcKernan96} \by S.~Keel and  J.~McKernan
    \paper Contractible extremal rays of $\overline{M}_{0,n}$
    \jour preprint alg-geom/9707016
    \yr 1996
\endref
\ref \key{Koll\'arMori98} \by J.~Koll\'ar and S.~Mori
     \book Birational geometry of algebraic varieties. With
the collaboration of C. H. Clemens and A. Corti
     \bookinfo Cambridge Tracts in Math. 
\vol 134
     \publ Cambridge University Press
\publaddr Cambridge
     \yr 1998
\endref
\ref \key{KontsevichManin96} \by M.~Kontsevich and Y.~Manin
     \paper Quantum cohomology of a product. With an
appendix by R. Kaufmann
     \jour  Invent. Math. 
     \vol 124 
     \yr 1996
     \pages 313--339
\endref
\ref \key{Logan00} \by A.~Logan
      \paper  Relations among divisors on the moduli space of 
curves with marked points
      \paperinfo preprint math.AG/0003104
      \yr 2000
\endref
\ref \key{LosevManin00} \by A.~Losev and Y.~Manin
      \paper New Moduli Spaces of Pointed Curves and Pencils
of Flat Connections
      \jour Michigna Math. J. 
      \vol 48 
      \yr 2000
      \pages 443--472
\endref
\ref \key{Moriwaki98} \by  A.~Moriwaki
     \paper Relative Bogomolov inequality and cone of semi-positive
divisors on $M_g$
     \jour J. Amer. Math. Soc.
     \vol 11
      \yr 1998
      \pages 569--600
\endref

\ref \key{Moriwaki00} \bysame
     \paper Nef divisors in codimension one on the moduli space of 
stable curves
     \paperinfo preprint math.AG/0005012
      \yr 2000
\endref

\ref\key{Moriwaki01} \bysame
     \paper The $Q$-Picard group of the moduli space of curves 
in positive characteristic
     \paperinfo preprint math.AG/0012016
     \yr 2000
\endref


\ref \key{Rulla00} \by W.~Rulla
      \book Birational Geometry of $\vmg 3.$
      \bookinfo Ph D. Thesis, Univ. of Texas at Austin
      \yr 2000
\endref
\endRefs
\enddocument